\input amstex
\magnification 1100
\NoBlackBoxes
\centerline { THE UNIRATIONALITY OF THE MODULI SPACES OF CURVES OF GENUS $\leq$ 14} \rm
\bigskip 
\centerline {ALESSANDRO VERRA}
\par
\centerline {Dipartimento di Matematica, Universit\`a Roma Tre}
 \bigskip \noindent
\bf 0. Introduction \rm  In this paper we prove that the moduli space of complex curves of genus 14 
is unirational. Our method applies more in general to the moduli space $\Cal M_g$ with $g \leq 14$,
so we use it to give new proofs of the unirationality of $\Cal M_g$ for $g = 11, 12, 13$. \par \noindent
The proof relies on linkage of curves in the projective space and on Mukai's description of canonical curves, of certain low genera, as
linear sections of a homogeneous space, ([M]). From these results of Mukai we deduce the unirationality of the Hilbert schemes of \it non
special, \rm smooth, irreducible curves of degree $d$ and genus $g \leq 10$ in $\bold P^r$, (see section 1). Then we use this property,
together with linkage, for proving our results. \par \noindent
\par \noindent To add some historical remarks we recall that the proof of the unirationality of $\Cal M_g$
goes back to Severi for $g \leq 10$, [Se]. The cases of genus $11,12,13$ were first proved  by Sernesi
and by Chang and Ran, [S1] for $g = 12$ and [CR1] for $g = 11,13$. Quite recently a proof which is in part
computational was given by Tonoli and Schreyer for $g = 11, 12, 13$, [TS]. \par \noindent
A conjecture of Harris and Morrison implies that $\Cal M_g$ has negative Kodaira dimension for $g \leq 22$. Our result implies such
a property for $\Cal M_{14}$, this was known up to now for $g \leq 13$ and $g = 15,16$, cfr. [HM2], [CR2],
[FP]. Of course things are different in higher genus: due to the fundamental results of Eisenbud, Harris and Mumford, $\Cal M_g$ has
non negative Kodaira dimension for $g \geq 23$ and it is of general type for $g \geq 24$. Recently Farkas showed that $\Cal M_{23}$ has
Kodaira dimension $\geq 1$, [F].
\par \noindent
Our starting point has been the following observation: fix a curve $D$ of
genus 14 with general moduli. On $D$ there are finitely many line bundles $L$ of degree 8 such that $h^0(L) = 2$. For 
each of them $\omega_D(-L)$ is very ample and defines an embedding
$$
D \subset \bold P^6.
$$
Now consider the vector space $V$ of quadratic forms vanishing on $D$: if $D$ is projectively normal then
$V$ has dimension 5, hence
$$
D \subset Q_1 \cap \dots \cap Q_5 
$$
where $Q_1 \dots Q_5$ are independent quadrics. If $Q_1 \dots Q_5$ define a complete intersection, then
$$
Q_1 \cap \dots Q_5 = C \cup D 
$$
where $C$ is a curve of degree 14. If $C$ is smooth and connected then its geometric genus is 8. In
section 5 we show the existence of a complete intersection of 5 quadrics
$$
C_o \cup D_o \subset \bold P^6 
$$
which satisfies all the previous assumptions. Using the general set up proved in sections 2, 3 and 4 we are also
able to deduce that $C_o$ is a non specially embedded, projectively normal curve and that the Petri map
$$
\mu: H^0(\omega_{D_o}(1)) \otimes H^0(\Cal O_{D_o}(1)) \to H^0(\omega_{D_o}) 
$$
is injective. Then, in their corresponding Hilbert schemes, $C_o$ and $D_o$ admit irreducible open neighborhoods
$\Cal C$ and $\Cal D$ parametrizing curves with the same properties. It follows from the results of section 1 that $\Cal C$ 
is unirational. On the other hand the injectivity of $\mu$ implies that the natural map $f: \Cal D \to \Cal M_{14}$
is dominant. On $\Cal C$ one can easily construct a Grassmann bundle
$$
\Cal G
$$
which is locally trivial in the Zariski topology and parametrizes pairs $(C,V)$ such that $C \in \Cal C$ and $V \subset H^0(\Cal I_C(2))$
is a 5-dimensional subspace. Since $\Cal C$ is unirational the same holds for $\Cal G$. Finally, the existence of the above  complete 
intersection $C_o \cup D_o$ makes possible to define a rational map
$$
\phi: \Cal G \to \Cal D 
$$
sending the general pair $(C,V)$ to $D$, where $C \cup D$ is the scheme defined by $V$ and $D$ is  a smooth, irreducible element
of $\Cal D$. It turns out that $\phi$ is birational, hence $f \cdot \phi: \Cal G \to \Cal M_{14}$ is dominant and $\Cal M_{14}$
is unirational. A more elaborated, from the technical point of view, version of this idea works for showing the unirationality of the
universal 5-symmetric product over $\Cal M_{12}$, of the universal 6-symmetric product over $\Cal M_{11}$ and finally of $\Cal M_{13}$. So
we are able to obtain in this way the unirationality of $\Cal M_g$, $g = 11, 12, 13$. \par \noindent
To continue with the example of genus 14 we give some more details on the way we prove the unirationality of $\Cal C$.
In the general case of the Hilbert scheme of non special curves of degree $d$ and genus $g \leq 10$ the
proof of the unirationality is analogous, (see section 2). It follows from Mukai's results that a general genus 8 canonical curve is a
linear section of the Pluecker embedding 
$$
G \subset \bold P^{14}
$$
of the Grassmannian of lines of $\bold P^5$, cfr. [M]. Assume $x = (x_1, \dots, x_8) \in G^8$ is general and let
$P_x$  be the space spanned by $x_1, \dots, x_8$. Then
$$
C_x =: P_x \cap G 
$$ 
is such a general canonical curve, moreover $C_x$ is endowed with the line bundle
$$
H_x \ =: \ \omega_{C_x}(x_1+ \dots+x_4-x_5-\dots -x_8) \in Pic^{14}(C_x).
$$
The pair $(C_x,H_x)$ defines a point in the universal Picard
variety $Pic_{14,8}$ and a rational map
$$
G^8 \to Pic_{14,8}. 
$$
It turns out that the latter is dominant, (see section 3), therefore $Pic_{14,8}$ is unirational. Then, with a little bit more
of effort, the unirationality of $\Cal C$ also follows. \par \noindent 
 \bigskip \noindent \it Some
frequently used notations and conventions.
\rm \par \noindent 
-$\Cal H_{d,g,r}$ denotes the \it restricted \rm Hilbert scheme. This is the
subscheme, in the Hilbert scheme of curves of
degree $d$ and genus $g$ in $\bold P^r$, parametrizing smoothable, connected, non degenerate curves. \par \noindent
- $W^r_d(D)$ is the Brill-Noether locus of all line bundles $H \in Pic^d(C)$ such that $h^0(H) \geq r+1$.
$\Cal W^r_{d,g}$ is the universal Brill-Noether locus over $\Cal M_g$. $Pic_{d,g}$ is the
universal Picard variety.
\par \noindent
-To simplify notations $X \cap Y$ will be the \it scheme theoretic \rm intersection of the schemes $X$
and $Y$, unless differently stated. \par \noindent
- $\Cal I_{X/Y}$ is the ideal sheaf of $X$ in $Y$. If $V$ is a vector space of sections of a line bundle then $\mid V \mid$
is the associated linear system. The dual of a vector bundle $\Cal E$ is $\Cal E^*$, $\bold P(\Cal E) =: Proj \ \Cal E^*$. 
\par \noindent - A nodal curve is a curve having ordinary nodes as its only singularities. 
\par \noindent (*) The present work has been supported by the
European research program EAGER and by the Italian research program GVA (Geometria delle variet\'a algebriche)
 \bigskip \noindent
\bf 1 Auxiliary unirationality results. \rm In this section we show some results of unirationality for the
Hilbert schemes $\Cal H_{d,g,r}$ in genus $g \leq 10$. More precisely let
$$
\Cal H^{ns}_{d,g,r} =: {\lbrace C \in \Cal H_{d,g,r} \ / \ \text {$C$ is smooth and $\Cal O_C(1)$ is non special }
\rbrace }, \tag 1.1
$$
where $3 \leq r \leq d-g$. As is well known the Zariski closure of $\Cal H^{ns}_{d,g,r}$ is the unique irreducible component of the Hilbert
scheme which dominates $\Cal M_g$. We will prove the following
\bigskip \noindent
(1.2) \bf THEOREM \it $\Cal H^{ns}_{d,g,r}$ is unirational for $7 \leq g \leq 10$.
\rm \bigskip \noindent
The extension of the theorem to genus $g \leq 6$ is easy: see 1.11. The theorem is an
application of the following description, due to Mukai,of a general canonical curve of genus $g = 7,8,9,10$ (cfr. [M]).
\bigskip \noindent
(1.3) \bf THEOREM \it For $g = 7,8,9,10$ there exists a rational homogeneous space $P_g  \subset \bold P^{dim P_g + g -
2}$ whose general curvilinear section is a general canonical curve.  
\rm \bigskip \noindent
Let $P_g$ be a homogeneous space as above, we consider the open subset
$$
U \subset P_g^g \tag 1.4
$$
of points $x = (x_1, \dots, x_g)$ such that $x_i \neq x_j (i \neq j)$ and moreover: \bigskip \noindent \it 
(i) the linear span $\bold P_x =: <x_1 \dots x_g>$ has dimension $g-1$, \par \noindent
(ii) $C_x =: \bold P_x \cap G$ is a smooth, irreducible canonical curve.
\rm \bigskip \noindent 
On $U$ we have a universal canonical curve 
$$
\pi: \Cal C \to U
$$
with fibre $C_x$ over $x$. $\Cal C$ contains the divisors $D_i = s_i(U)$, where $s_i: U \to \Cal C$ is the section sending $x$ to
$x_i$. For any $d \geq g+3$ we fix \it non zero \rm integers \rm $n_1, \dots, n_g$ such that 
$$
d = n_1 + \dots + n_g + 2 - 2g. 
$$
Then we consider the sheaf
$$
\Cal H =: \omega_{\pi}(n_1D_1+ \dots + n_gD_g), \tag 1.5
$$
where $\omega_{\pi}$ denotes the relative cotangent sheaf of $\pi$. For any $x \in U$ we have
$$
\Cal H \otimes \Cal O_{C_x}Ê= \omega_{C_x}(n_1x_1+ \dots +n_gx_g).
$$
Let $C$ be a curve of genus $g$, then the Abel map $a: C^g \to Pic^g(C)$ is surjective. As is well known 
this property generalizes as follows: \bigskip \noindent
(1.6) \bf LEMMA \it Let $n_1, \dots, n_g$ be non zero integers, then the map
$$
a_{n_1 \dots n_g}: C^g \to Pic^d(C)
$$
sending $(x_1, \dots, x_g)$ to $\omega_C(n_1x_1+ \dots + n_gx_g)$ is surjective. \rm \bigskip \noindent
(1.7) \bf LEMMA \it Let $x$ be general in $U$, then $\Cal H \otimes \Cal O_{C_x}$ is non
special and very ample. \rm \bigskip \noindent
PROOF Let $C = C_x$, by Mukai's result we can assume that $C$ has general moduli. By the previous lemma, each
$L \in Pic^d(C)$ is isomorphic to $\omega_C(\Sigma n_iy_i)$, for some $y = (y_1, \dots, y_g) \in C^g$.
Therefore, keeping $C$ fixed and possibly replacing $x$ by $y$, we can assume that $\Cal H \otimes \Cal O_C$ is general 
in $Pic^d(C)$. On a general curve of genus $g$ a general line bundle of degree $d \geq g+3$ is very ample, as follows
from [ACGH] Theorem 1.8 p.216. Moreover such a line bundle is also non special.
\bigskip \noindent 
Finally, up to shrinking the open set $U$, we can assume that: \it \bigskip \noindent
(iii) $Proj (\pi_* \Cal H) = U \times \bold P^{d-g}$, \par \noindent
(iv) $\Cal C$ is embedded in $U \times \bold P^{d-g}$ by the map associated to the tautological sheaf, \par \noindent
(v) $\Cal H \otimes \Cal O_{C_x}$ is non special, $\forall x \in U$. \rm \bigskip \noindent
(1.8) \bf LEMMA \it The natural morphism $r: U \to Pic_{d,g}$ is dominant. \rm \bigskip \noindent
PROOF Let $x \in U$, by definition $r(x)$ is the isomorphism class of $(C_x, \Cal H \otimes \Cal O_{C_x})$ in the
universal Picard variety $Pic_{d,g}$.  To prove that $r$ is dominant, it suffives to show that $r(U)$ intersects a general fibre of
the forgetful map $f: Pic_{d,g} \to \Cal M_g$ along a dense subset. A general fibre of $f$ is $Pic^d(C)$, where $C$ is general
of genus $g$. By Mukai's theorem $C$ is biregular to a section $C_o = \bold P_o \cap P_g$,
for some point $o = (o_1, \dots, o_g)$ in $C_o^g \cap U$. Note that
$r/C_o^g \cap U: C_o^g \to Pic^d(C)$ associates to $x = (x_1, \dots, x_g)$ the line bundle $\Cal H \otimes \Cal O_{C_x}Ê= \Cal
O_C(n_1x_1+ \dots + n_gx_g)$. Then such a map extends to the surjection $a_{n_1 \dots n_g}$ considered in lemma 2.6 and hence
$r(C_o^g \cap U)$ is dense. \bigskip \noindent
Let
$$
p: \Cal H^{ns}_{d,g,d-g} \to Pic_{d,g} \tag 1.9
$$
be the morphism sending $C$ to the isomorphism class of the pair $(C,\Cal O_C(1))$ and let
$$
q: U \times PGL(d-g+1) \to \Cal H^{ns}_{d,g,d-g} \tag 1.10
$$
be the morphism sending $(x,\alpha)$ to $\alpha(C_x)$, observe that
$$
p \cdot q (U \times \lbrace id \rbrace) = r(U).
$$
Therefore, by lemma 1.8, $\ p \cdot q: U \times PGL(d-g+1) \to Pic_{d,g}$ is dominant. 
\bigskip \noindent
PROOF OF THEOREM 1.2 We first show the case $r = d-g$. We know from 1.8 and 1.10 that both $p(q(U \times PGL(r+1)))$ and $p(\Cal
H^{ns}_{d,g,d-g})$ contain a non empty open set of $Pic_{d,g}$. Therefore
$$
A =: p^{-1}(p((U \times PGL(d-g+1)) \cap p(\Cal H^{ns}_{d,g,d-g}))
$$
contains a non empty open set. Then, since $\Cal H^{ns}_{d,g,r}$ is irreducible, it suffices to show that
$$
A \subset q(U \times PGL(r+1)).
$$
Indeed this implies that $q$ is dominant and hence that $\Cal H^{ns}_{d,g,r}$ is unirational. To prove the above inclusion we observe 
that: $C \in A$ $\Longrightarrow$ $p(C)$ $\in$ $p(q(U \times PGL(d-g+1)))$ $\Longrightarrow$ $C$ is smooth and $\Cal O_C(1)$ is non 
special $\Longrightarrow$ $C$ is linearly normal $\Longrightarrow$ $p^{-1}(p(C)) = P_C$, where $P_C$ is the $PGL(d-g+1)$-orbit of $C$. 
But then there  exists $\beta \in PGL(r+1)$ such that $\beta(C) = q(x,\alpha)$, that is $C = q(x, \beta^{-1} \alpha)$. Let us complete 
the proof with the case $r < d-g$: consider the space $M_r$ of all linear projections  $\bold P^{d-g} \to \bold P^r$ and the map 
$s: Y \times M_r \to \Cal H^{ns}_{d,g,r}$ sending  $(C,\alpha)$ to $\alpha(C)$. It is easy to see that $s$ is dominant, therefore $\Cal
H^{ns}_{d,g,r}$ is unirational too.
\bigskip \noindent
(1.11) \bf REMARK \rm Let $4 \leq g \leq 6$ and let $P_g \subset \bold P^{g+1}$ be a fixed, general Fano threefold of index one 
and genus $g$. A direct count of parameters shows that a general canonical curve of genus $g$ is a curvilinear section of $P_g$.
Using this property one can show, with exactly the same proof as above, that $\Cal H^{ns}_{d,g,r}$ is unirational if $4 \leq g \leq 6$.
We leave to the reader the extension of the result to the case $g \leq 3$.
\bigskip \noindent
\bf 2. General set up. \rm In this section we build up a somehow general strategy to apply the previous unirationality results. The
basic idea is to consider families of irreducible curves $D$ of genus $g'$ which are linked to a general $C \in \Cal H^{ns}_{d,g,r}$ 
by a complete intersection of fixed type $(f_1, \dots, f_{r-1})$. These families are unirational. In some cases they dominate $\Cal
M_{g'}$. We start with the following \bigskip \noindent 
(2.1) \bf DEFINITION \it Let $3 \leq r \leq d-g$ and let $g \leq 10$, then
$$
\Cal C_{d,g,r} \ =: \ \lbrace \ C \in \Cal H^{ns}_{d,g,r} \ / \ \rho_f \ \text {has maximal rank for each $f$} \ \rbrace
$$
where $\rho_f: H^0(\Cal I_{C/\bold P^r}(f)) \to H^0(\Cal O_C(f))$ is the restriction map. \rm \bigskip \noindent
By semicontinuity $\Cal C_{d,g,r}$ is open, it is non empty because the maximal rank condition is generically satisfied 
in the cases we are considering. Due to the results of the previous section $\Cal C_{d,g,r}$ is irreducible and unirational. We fix a
sequence of integers
$$
\sigma = (f_1, \dots, f_s, k_1, \dots, k_s) \tag 2.2
$$
satisfying 
$$
1 \leq \ k_i  \leq  n_i \ \ \text {and} \ \ k_1 + \dots + k_s = r-1
$$ 
where $n_i$ is the constant value of $h^0(\Cal I_{C/\bold P^r}(f_i))$ when $C$ moves in $\Cal C_{d,g,r}$. Then we consider the Ideal 
sheaf $\Cal J$ of the universal curve
$$
\Cal C \subset \Cal C^{ns}_{d,g,r} \times \bold P^r
$$
and the natural projections  $p_1: \Cal C \to \Cal C_{d,g,r}$ and $p_2: \Cal C \to \bold P^r$. By Grauert theorem the sheaf
$$
\Cal F_i =: p_{1*} (\Cal J \otimes  p_2^* \Cal O_{\bold P^r}(f_i)) \tag 2.3
$$
is a vector bundle on $\Cal C_{d,g,r}$, the fibre at $C$ is the space $H^0(\Cal I_{C/\bold P^r}(f_i))$. For each $k \geq 1$ we can
also consider the Grassmann bundle
$$
u_k: G(k,\Cal F_i) \to \Cal C^{ns}_{d,g,r} \tag 2.4
$$
defined by $k$ and $\Cal F_i$. 
\bigskip \noindent
(2.5) \bf DEFINITION \it $\Cal G^{\sigma}_{d,g,r}$ is the fibre product $G(k_1,\Cal F_1) \times_{u_{k_1}} \dots \times_{u_{k_s}} G(k_s,
\Cal F_s)$ over $\Cal C_{d,g,r}$. \rm \bigskip \noindent
$\Cal G^{\sigma}_{d,g,r}$ is birational to the product $\Cal C_{d,g,r} \times G(k_1,n_1) \times \dots \times G_(k_s,n_s)$, therefore the
next statement is immediate.  
\bigskip \noindent 
(2.6) \bf PROPOSITION \it $\Cal G^{\sigma}_{d,g,r}$ is irreducible and unirational. \rm \bigskip \noindent
\bigskip \noindent 
A point of $\Cal G^{\sigma}_{d,g,r}$ is a sequence $(C,V_1, \dots, V_s)$, where $C \in \Cal C^{ns}_{d,g,r}$ and $V_i$ is a
$k_i$-dimensional subspace of $H^0(\Cal I_{C/\bold P^r}(f_i))$. Let
$$
B \subset \bold P^r
$$
be the scheme defined by the set of homogeneous forms $V_1 \cup \dots V_s$. Since $dim \ V_1 + \dots + dim \ V_s = r-1$
it is possible that $B$ is a curve: in this case $B$ is a complete intersection. 
\bigskip \noindent
(2.7) \bf DEFINITION \it (1) A point $(C,V_1, \dots, V_s) \in \Cal G^{\sigma}_{d,g,r}$ is a key-point if 
$$
B = C \cup D
$$
is a nodal curve and the component $D$ is smooth, irreducible, non degenerate. \bigskip \noindent
(2) $\Cal G^{\sigma}_{d,g,r}$ satisfies the key condition if the set of its key-points is non empty.
\rm \bigskip \noindent
Clearly the set of the key-points of $\Cal G^{\sigma}_{d,g,r}$ is open. Let $C \cup D$ be a nodal complete intersection as above,
from now on we will keep the following notations: 
$$
d' = deg(D) \ , \ g' = p_a(D) \ , \ \text { $n$ $=$ cardinality of $Sing \ B$. } \tag 2.8
$$
The numbers $d'$, $g'$ and $n$ can be readily computed from $(d,g,\sigma)$, we have: 
$$
d + d' = f_1^{k_1} \cdots f_{s}^{k_s}, 
$$
$$
(g-g') = \frac 12 (k_1f_1+ \dots +k_sf_s -r-1)(d-d')
$$
and
$$
n = (k_1f_1+ \dots +k_sf_s -r-1)d + 2 - 2g, 
$$
see [Fu] p.159 example 9.1.12. Finally we fix the notations for some natural maps: \bigskip \noindent
(2.9) \bf DEFINITION \it  Assume that $\Cal G^{\sigma}_{d,g,r}$ satisfies the key condition, then 
$$
\gamma_{d,g,r}: \Cal G^{\sigma}_{d,g,r} \to \Cal H_{d',g',r}
$$
is the map sending a point $(C,V_1, \dots,V_s)$ as above to $D$. \rm \bigskip \noindent
The natural map from $\Cal H_{d,g,r}$ to the universal Brill-Noether locus will be denoted as
$$
\beta_{d,g,r}: \Cal H_{d,g,r} \to \Cal W^r_{d,g} \tag 2.13
$$
and the forgetful map from $\Cal W^r_{d,g}$ to the moduli space will be denoted as 
$$
\alpha_{d,g,r}: \Cal W^r_{d,g} \to \Cal M_{g} \tag 2.14.
$$
The next statement only summarizes our program for showing the unirationality of some moduli spaces, the proof is immediate. 
\bigskip \noindent
(2.10) \bf PROPOSITION \it Assume $\Cal G^{\sigma}_{d,g,r}$ satisfies the key condition and that the image of
$$
\gamma_{d,g,r}: \Cal G^{\sigma}_{d,g,r} \to \Cal H_{d',g',r}
$$
dominates $\Cal M_{g'}$, then $\Cal M_{g'}$ is unirational. \rm \bigskip \noindent 
\bf 3. Some useful criteria I.  \rm In this section we prove sufficent conditions to ensure that the key condition holds for 
$\Cal G^{\sigma}_{d,g,r}$. This criterion is an elementary version of more general known properties, we show it for completeness and some
lack of reference. We recall that a subvariety
$$
Y \subset \bold P^r
$$
is a scheme theoretic intersection of hypersurfaces of degree $f$ iff $\Cal I_{Y/\bold P^r}(f)$ is globally generated. Throughout
all the section we will assume that
$$
Y =: Y^c \cup Z \tag 3.1
$$
where: \par \noindent
- $Y^c$ is an equidimensional variety of codimension $c$, which is locally complete intersection with at most
finitely many singular points. \par \noindent
- $Z$ is disjoint from $Y^c$ and it is either smooth 0-dimensional or empty. \bigskip \noindent
(3.2) \bf PROPOSITION \it Assume that $\Cal I_{Y/\bold P^r}(f)$ is globally generated and that $dim \ Y \leq 3$. Then there exists a
complete intersection of $c$ hypersurfaces $Q_1 \dots Q_c$ of degree $f$ such that: \par \noindent
- either $Q_1 \cap \dots \cap Q_c = Y$ \par \noindent
- or $Q_1 \cap \dots \cap Q_c = X \cup Y$ and moreover: (1) $X$ is smooth and contains $Z$, (2) $X \cap Y^c$ is smooth and equidimensional
of codimension $c+1$. 
\rm \bigskip \noindent  PROOF Let $I =: H^0(\Cal I_{Y/\bold P^r}(f))$, we denote by $G^c$ the Grassmannian of codimension $c$ subspaces of
$I$ and by $G_c$ the Grassmannian of $c$ dimensional subspaces. We assume that $Y$ is not a complete intersection of $c$ hypersurfaces
of degree $f$: otherwise there is nothing else to show. Let
$$
\sigma: P \to \bold P^r
$$
be the blowing up of $Y^c$. Then the strict transform of $\mid I \mid$ by $\sigma$ is $\mid fH - E \mid$, where $E$ is the exceptional
divisor of $\sigma$ and $H$ is the pull-back of a hyperplane. Since $\Cal I_{Y/\bold P^r}(f)$ is globally generated and $Y = Y^c \cup Z$,
the base locus of $\mid fH - E \mid$ is $\sigma^{-1}(Z)$. Let $V \in G_c$ be general and let
$$
B_V
$$
be the base locus of the strict transform of $\mid V \mid$ on $P$. By Bertini theorem we can assume that $B_V$ is smooth: this
follows because $\sigma^{-1}(Z)$ is smooth and finite. Moreover we can assume that $B_V$ intersects transversally the exceptional divisor
$E$. Since $Y^c$ is locally complete intersection, $E$ is a projective bundle with fibre of dimension $c-1$. Then, since $V$ is
general of codimension $c$ and $Sing \ Y^c$ is finite, we can also assume that $B_V \cap \sigma^{-1}(Sing \ Y^c)$  is empty. We claim that
$$
\sigma/B_V: B_V \to \bold P^r
$$
is an embedding: this is obvious on $B_V-(B_V \cap E)$. To complete the proof consider $p \in B_V \cap E$ and $F =$
 $\sigma^{-1}(\sigma(p))$. $B_V$ is the complete intersection of $c$ independent divisors $D_1, \dots D_c$  of 
$\mid fH - E \mid$. Since $\Cal O_F(fH-E) \cong \Cal O_{\bold P^{c-1}}(1)$ the intersection scheme $B_V \cap F$ is a linear  space. 
This must be 0-dimensional because $B_V$ is transversal to $E$, hence $\sigma/B_V$ is an embedding at $p$. Let
$$
X_V = \sigma(B_V)
$$
then $X_V$ is smooth and moreover we have that
$$
X_V \cup Y
$$
is complete intersection of the $c$ hypersurfaces $\sigma(D_1), \dots, \sigma(D_c) \in \mid \Cal I_{Y/\bold P^r}(f) \mid$. \par \noindent
Finally we show that $X_V \cap Y^c$ is smooth. Let $y \in Y^c - Sing \ Y$ and let
$$
I_y = \lbrace q \in I \ / \ q \in m^2_y \rbrace.
$$
Since $Y$ is scheme theoretic intersection of hypersurfaces of degree $f$, the space $I_y$ has codimension $c$ in $I$. This defines
a morphism
$$
\phi: Y^c - Sing \ Y \to G^c
$$
sending $y$ to $I_y$. For any $V \in G_c$ we can consider the Schubert cycle
$$
\sigma_V = \lbrace L \in G^c \ / \ dim (L \cap V) \geq 1 \rbrace.
$$
It is well known that the singular locus of $\sigma_V$ is
$$
Sing \ \sigma_V = \lbrace L \in G^c \ / \ dim \ (L \cap V) \geq 2 \rbrace
$$
and moreover that $Sing \ \sigma_V$ has codimension 4 in $G_c$. Since $dim \ Y \leq 3$ we have
$$
Sing \ \sigma_V \cap \phi(Y^c - Sing \ Y) = \emptyset
$$
for a general $V$. Then, by the transversality of a general $\sigma_V$, we can assume that $\phi^{-1}(\sigma_V)$ is 
smooth of codimension $c+1$. On the other hand it turns out that
$$
\phi^{-1}(\sigma_V) = Y^c \cap X_V.
$$
Indeed: $y \in Y^c \cap X_V$ $\Leftrightarrow$ $y \in Y^c-Sing \ Y$ and $B_V$ is singular at $y$ $\Leftrightarrow$
$y \in Y^c - Sing \ Y$ and $dim \ I_y \cap V = 1$.
 \bigskip \noindent
In particular we will apply the lemma when $Y^c$ is a nodal curve. So we point out the following
\bigskip \noindent
(3.3) \bf PROPOSITION \it Let $Y = C \cup Z$, where $C$ is a nodal curve and $Z$ is a smooth, 0-dimensional scheme disjoint
from $C$. Assume that $\Cal I_{Y/\bold P^r}(f)$ is globally generated, then: 
\par \noindent 
(1) $Y$ lies in a smooth surface $S$ which is complete intersection of $r-2$ hypersurfaces of degree $f$. 
\par \noindent 
(2) $\mid fH - C \mid$ is base-point-free, so that a general $D \in \mid fH - C \mid$ is transversal to $C$. $D$ is connected if
$D^2 > 0$. ($H$ $=$ hyperplane section of $S$).
\rm \bigskip \noindent
PROOF By proposition 3.2 there exists a nodal complete intersection of $r-1$ hypersurfaces of degree $f$
$$
C \cup D 
$$
such that $D$ is smooth and contains $Z$. Let $V = H^0(\Cal I_{C \cup D/ \bold P^r}(f))$ and let  $x \in C \cup D$: if $x$ is smooth 
no $Q \in \mid V \mid$ is singular at $x$. If $x$ is singular then $x$ is a node and there exists exactly one $Q$ singular at $x$.
Since $Sing C \cup D$ is finite, the base locus of a general hyperplane in $\mid V \mid$ is a smooth surface $S$ as required. This
shows (1), (2) follows immediately. \bigskip \noindent
(3.4) \bf PROPOSITION \it Let $Y = C \cup Z$ be as in \rm 3.3 and let $\sigma = (f,r-1)$. \it Assume that: \par \noindent
(1) $\Cal I_{Y/\bold P^r}(f)$ is globally generated and $C \in \Cal C_{d,g,r}$. \par \noindent
(2) $d \leq f^{r-1} - f^{r-2}$ and $d'(f^{r-2}-r-1) < 2g'-2$.  \par \noindent
Then $\Cal G^{\sigma}_{d,g,r}$ satisfies the key condition. \rm \bigskip
\noindent PROOF By 3.3 $Y \subset S$, where $S$ is a smooth complete intersection of $r-2$ hypersurfaces of degree $f$. Moreover
there exists a smooth $D \in \ \mid fH - C \mid$ transversal to $C$. Conditions in (2) are equivalent to $(H-D)H \leq 0$ and $D^2
> 0$. Hence $D$ is connected, not degenerate and the statement follows.  \bigskip
\noindent
\bf 4. Some useful criteria II. \rm Now we want to give some sufficient conditions, on some of the unirational Grassmann bundles $\Cal
G^{\sigma}_{d,g,r}$, so that 
$\Cal G^{\sigma}_{d,g,r}$ dominates the moduli space $\Cal M_{g'}$. \bigskip \noindent
(4.1) \bf LEMMA \it Let $C \cup D \subset \bold P^r$ be a nodal complete intersection of $r-1$
hypersurfaces of degree
$f =
\frac {r+2}{r-2}$. Then:  \par \noindent  
(1) $C$ is $f$-normal iff $D$ is linearly normal and $D$ is $f$-normal iff $C$ is linearly normal. \par \noindent 
(2) $\Cal O_C(1)$ is non special iff $h^0(\Cal I_{D/\bold P^r}(f)) = r-1$. \par \noindent
(3) $C$ is non degenerate iff $\Cal O_D(f)$ is non special and $D$ is non degenerate iff $\Cal O_C(f)$ is non special.
\rm \bigskip \noindent PROOF As in 3.3 $C \cup D$ is contained in a smooth complete intersection $S$ of $r-2$ hypersurfaces
of degree $f$. \par \noindent
(1) Let $H$ be a hyperplane section of $S$, the assumption $f = \frac {r+2}{r-2}$ simply means that $H$ is a canonical divisor. It
follows from the standard exact sequence
$$
0 \to \Cal O_S(H-D) \to \Cal O_S(H) \to \Cal O_D(H) \to 0
$$ 
that $D$ is linearly normal iff $h^1(\Cal O_S(H-D)) = 0$. Since $H$ is canonical we have $h^1(\Cal O_S(H-D))$ $=$ $h^1(\Cal
O_S(D)) = h^1(\Cal O_S(fH-C))$. Hence $D$ is linearly normal iff $C$ is $f$-normal. Since $C+D \sim fH$, the second equivalence follows
exchanging $D$ with $C$. \par \noindent  
(2)  At first we remark that $h^0(\Cal I_{D/\bold P^6}(f)) = h^0(\Cal O_S(fH-D)) + r-2 =
h^0(\Cal O_S(C)) + r-2$. Secondly the standard exact sequence
$$
0 \to \Cal O_S(H-C) \to \Cal O_S(H) \to \Cal O_C(1) \to 0
$$
yelds $h^1(\Cal O_C(1)) = h^2(\Cal O_S(H-C)) - h^2(\Cal O_S(H)) = h^0(\Cal O_S(C)) - 1$. Hence: $h^0(\Cal I_{D/\bold P^6}(2)) =
r-1$ $\Leftrightarrow$ $h^0(\Cal O_S(C)) = 1$ $\Leftrightarrow$ $h^1(\Cal O_C(1)) = 0$. \par \noindent
(3) Consider the standard exact sequence
$$
0 \to \Cal O_S(C) \to \Cal O_S(fH) \to \Cal O_D(f) \to 0
$$
and its associated long exact sequence. Since $f \geq 1$, we have $h^1(\Cal O_S(fH)) = h^2(\Cal O_S(fH)) = 0$ and hence $h^1(\Cal O_D(f))$
$=$ $h^2(\Cal O_S(C))$ $=$ $h^0(\Cal O_S(H-C))$. This implies the first equivalence. Again the second one follows with by exchanging
$C$ and $D$.
 \bigskip \noindent  
(4.2) \bf PROPOSITION \it  Let $\sigma = (f,r-1)$, $f = \frac {r+2}{r-2}$ and $r = d-g$. If $\Cal G^{\sigma}_{d,g,r}$ satisfies
the key assumption. Then the image of the map $\ \ \beta_{d',g',r} \cdot \gamma_{d,g,r}: \Cal G^{\sigma}_{d,g,r} \to \Cal W^r_{d',g'} \  $
is open.
\rm \bigskip \noindent
PROOF Under the assumption a general $C \in \Cal C_{d,g,r}$ is linked to a smooth, irreducible, non degenerate curve $D$ by a complete
intersection of $r-1$ hypersurfaces of degree $f$. We remark that $C$ is both $f$-normal and linearly normal. This follows because the
restriction map
$$
\rho_m: H^0(\Cal O_{\bold P^r}(m)) \to H^0(\Cal O_C(m))
$$
has maximal rank and $h^0(\Cal I_{C/\bold P^4}(f)) > 0$. Hence $\rho(f)$ is surjective and $C$ is $f$-normal. On the other hand 
$r = d - g$ implies $h^0(\Cal O_C(1)) = r+1$, because $\Cal O_C(1)$ is non special. Then $\rho(1)$  is an isomorphism and $C$ is 
linearly normal. $D$ has the following properties:  \par \noindent (i)
$D$ is linearly normal, (ii) $D$ is $f$-normal, (iii) $h^0(\Cal I_{D/\bold P^r}(f)) = r-1$, (iv) $\Cal O_D(f)$ is non special. 
\par \noindent
This follows from the previous lemma 4.1. Using the same lemma it is easy to see that the set 
$$
U' = \lbrace D' \in \Cal H_{d',g',r} \ / \ \text {$D'$ is smooth, irreducible, non degenerate and satisfies (i), \dots, (iv) } \rbrace,
$$
is open. Let $D' \in U'$, assume that: (v) the scheme defined by $V' =: H^0(\Cal I_{D/\bold P^r}(f))$ is a complete intersection $C' \cup
D'$, where $C'$ is smooth, irreducible. Then, by lemma 4.1, $C' \in \Cal C_{d,g,r}$ and hence $D' = \gamma_{d,g,r}(C',V')$. Conversely,
if $D' = \gamma_{d,g,r}(C',V')$, then $D'$ satisfies (v). Thus condition (v) defines an open set $U \subset U'$ and $U$ is
the image of $\gamma_{d,g,r}$. $U$ is invariant under the action of $PGL(r+1)$, moreover each $D \in U$ is linearly normal. This
implies that $\beta_{d',g',r}(U) = U/PGL(r+1)$ and that $\beta_{d',g',r}(U)$ is open in $\Cal W^r_{d',g'}$.
\bigskip \noindent 
However recall that $\Cal W^r_{d',g'}$ could be reducible, even if the Brill-Noether number $\rho(d',g',r)$ is $\geq 0$. So it could
happen that $\gamma_{d,g,r} \cdot \beta_{d',g',r}$ is not dominant and that its image does not dominate $\Cal M_{g'}$. 
\bigskip \noindent
(4.3) \bf DEFINITION \it Let $x \in \Cal W^r_{d',g'}$ be the moduli point of the pair $(D,L)$, $L \in Pic^{d'}(D)$. We will say
that $x$ is Petri general if the Petri map $\mu: H^0(\omega_D(-L)) \otimes H^0(L)) \to H^0(\omega_D)$ is injective. The main universal
Brill-Noether locus is the open set
$$
\Cal U^r_{d',g'} = \lbrace \ x \in \Cal W^{r}_{d',g'}Ê\ / \ \text {$x$ is Petri general} \ \rbrace.
$$
\rm \bigskip \noindent
Let  $\rho(d',g',r) \geq 0$, by the main theorems of the Brill-Noether theory $\Cal U^r_{d',g'}$ is irreducible  and dominates $\Cal M_g$ 
via the natural map. This motivates the previous definition. 
\bigskip \noindent (4.4) \bf LEMMA \it Let $C \cup D \subset \bold P^r$ be a nodal complete intersection of $r-1$
hypersurfaces of degree $f =
\frac {r+2}{r-2}$ and let $r = d-g$. Assume that $C \in \Cal C_{d,g,r}$ then the multiplication
$$
\mu_C: H^0(\Cal O_{\bold P^r}(1)) \otimes H^0(\Cal I_{C/\bold P^r}(f)) \to H^0(\Cal I_{C/\bold P^r}(f+1))
$$
has maximal rank iff it has maximal rank the Petri map
$$
\mu: H^0(\omega_D(-1)) \otimes H^0(\Cal O_D(1)) \to H^0(\omega_D).
$$
In particular $\mu$ has maximal rank if the ideal of $C$ is generated by forms of degree $f$.
\rm\bigskip \noindent
PROOF We already know that $C \cup D \subset S$, where $S$ is a smooth canonical surface which is a complete intersection of
$r-2$ hypersurfaces of degree $f$. Then it holds $\omega_D(-1) \cong \Cal O_D(D)$. Moreover the standard exact sequence
$$
0 \to \Cal O_S \otimes H^0(\Cal O_S(H)) \to \Cal O_S(D) \otimes H^0(\Cal O_S(H)) \to \Cal O_D(D) \otimes H^0(\Cal O_S(H)) \to 0
$$
induces the exact commutative diagram
$$ 
\CD
0 @>>> {H^0(\Cal O_S(H))} @>>> {H^0(\Cal O_S(D)) \otimes H^0(\Cal O_S(H))} @>>> {H^0(\omega_D(-1)) \otimes H^0(\Cal O_D(1))} @>>> 0 \\
@. @VV{id}V @VV{\mu_S}V @VV{\mu}V \\
0 @>>> {H^0(\Cal O_S(H))} @>>> {H^0(\Cal O_S(H+D))} @>>> {H^0(\omega_D}) @>>> 0 \\
\endCD
$$
Since the left vertical arrow is an isomorphism we have $rk \ \mu_S = rk \ \mu + r+1$, hence $rk \ \mu_S$ is maximal iff
$rk \ \mu$ is maximal. Now observe that $\Cal I_{C/S}(f) \cong \Cal O_S(D)$ and consider the exact sequence
$$
0 \to V \otimes \Cal I_{S/\bold P^r}(f) \to V \otimes \Cal I_{C/\bold P^r}(f) \to V
\otimes 
\Cal O_S(D)
\to 0,
$$
where $V = H^0(\Cal O_{\bold P^r}(1))$. The sequence induces the exact diagram
$$ 
\CD
0 @>>> {V \otimes H^0(\Cal I_{S/\bold P^r}(f)} @>>> {V \otimes H^0(\Cal I_{C/\bold
P^r}(f)} @>>> {V \otimes H^0(\Cal O_S(D))} @>>> 0
\\ @. @VVV @VV{\mu_C}V @VV{\mu_S}V \\
0 @>>> {H^0(\Cal I_{S/\bold P^r}
(f+1))} @>>> {H^0(\Cal I_{C/\bold P^r}(f+1))} @>>> {H^0(\Cal O_S(H+D)}) @>>> 0 \\
\endCD
$$
Since $S$ is a complete intersection of hypersurfaces of degree $f$, the left vertical arrow is an isomorphism. Then $\mu_S$
has maximal rank iff $\mu_C$ has maximal rank and the first statement follows. The second statement is an obvious consequence: 
if the homogeneous ideal of $C$ is generated by forms of degree $f$ then $\mu_C$ is surjective, hence $\mu$ has maximal rank.  
\bigskip \noindent
(4.5) \bf THEOREM \it Let $\sigma = (f,r-1)$, $f = \frac {r+2}{r-2}$ and $r = d-g$. Assume $\rho(d',g',r) \geq 0$ and that
one of the following conditions hold: \par \noindent
(1) The ideal of a general $C \in \Cal C_{d,g,r}$ is generated by forms of degree $f$ and $2g'-2 > d' > f^{r-2}$. \par
\noindent (2) $\Cal G^{\sigma}_{d,g,r}$ satisfies the key condition and $\Cal W^r_{d',g'}$ is irreducible. \par \noindent
Then both $\Cal U^r_{d',g',r}$ and $\Cal M_{g'}$ are unirational.
\rm \bigskip \noindent   
PROOF Assume condition (1) holds and consider a general $(C,V) \in \Cal G^{\sigma}_{d,g,r}$. Then the scheme defined by $V$ is a
nodal complete intersection $C \cup D$, where $D$ is a smooth curve of arithmetic genus $g'$ and degree $d'$. It is easy to show
that the assumption $2g'-2$ $>$ $d' > f^{r-2}$ implies that $D$ is irreducible and non degenerate. Hence
$\Cal G^{\sigma}_{d,f,r}$ satisfies the key condition.  Notice also that
$$
\rho(d',g',r) \geq 0 \Longrightarrow h^0(\Cal O_D(1))h^0(\omega_D(-1)) \leq g',
$$
otherwise the Petri map $\mu$ would never be injective. Condition (1) implies that the multiplication 
$$
\mu_C: H^0(\Cal O_{\bold P^r}(1)) \otimes H^0(\Cal I_{C/\bold P^r}(f)) \to H^0(\Cal I_{C/\bold P^r}(f+1))
$$
is surjective. Then, by lemma 4.4, $\mu$ has maximal rank and hence it is injective. Therefore $\beta_{d',g',r}(D)$ is a point of the 
main universal Brill-Noether locus $\Cal U^r_{d',g'}$ and the image of
$$
\beta_{d',g',r} \cdot \gamma_{d,g,r}:  \Cal G^{\sigma}_{d,g,r} \to \Cal W^r_{d',g'}
$$
is contained in $\Cal U^r_{d',g'}$. On the other hand it follows from proposition 4.2 that such a image is open. Since $\Cal U^r_{d',g'}$
is irreducible and dominates $\Cal M_{g'}$ the statement follows. \par \noindent
Finally assume that (2) holds. Then the image of the above map is open and also dense in $\Cal W^r_{d',g'}$. Moreover $\Cal W^r_{d',g'}$
dominates $\Cal M_{g'}$ because $\rho(d',g',r) \geq 0$. Hence the statement follows again. \bigskip \noindent
\bf 5 Curves of degree 14 and genus 8 in $\bold P^6$. \rm Now we want to prove that the homogeneous ideal of a general curve
$$
C \in \Cal C_{14,8,6}
$$
is generated by quadrics. We start with a smooth, non degenerate Del Pezzo surface 
$$
Y \subset \bold P^6 \tag 5.1
$$ 
of degree 6. It is well known that $\Cal I_{Y/\bold P^6}(2)$ is globally generated. Then, by lemma 3.5, there exists a reducible,
nodal complete intersection of 4 quadrics
$$
X \cup Y \tag 5.2
$$
where $X$ is a smooth, irreducible surface of degree 10 and 
$$
F = X \cap Y \tag 5.3
$$
is a smooth curve. We have 
$$
\Cal O_{X \cup Y}(1) \cong \omega_{X \cup Y}
$$
for the dualizing sheaf of $X \cup Y$, moreover it holds
$$
\omega_X(F) \cong \omega_{X \cup Y} \otimes \Cal O_X \ , \ \omega_Y(F) \cong \omega_{X \cup Y} \otimes \Cal O_Y. \tag 5.4
$$
Since $Y$ is a Del Pezzo, it follows that $F \in \mid \Cal O_Y(2) \mid$ is a quadratic section of $Y$. In particular $F$ is a smooth
canonical curve of genus $7$ in $\bold P^6$. We point out that $F$ is non trigonal. This follows because a trigonal canonical curve has
infinitely many trisecant lines. Since $Y$ is intersection of quadrics, these lines would be contained in $Y$: a clear contradiction. 
Let $H_X \cup H_Y$ be a general hyperplane section of $X \cup Y$, with $H_X \in \mid \Cal O_X(1) \mid$ and $H_Y \in
\mid \Cal O_Y(1) \mid$. Since $F$ has degree 12, it follows
$$
17 = p_a(H_X \cup H_Y) = p_a(H_X) + p_a(H_Y) + 11 
$$
and hence $p_a(H_X) = 5$. \bigskip \noindent
(5.5) \bf PROPOSITION \it $X$ is rational and projectively normal. Moreover the homogeneous ideal of $X$ is generated by quadrics.
\rm \bigskip \noindent
PROOF To prove that a surface is projectively normal it suffices to show the same property for a smooth, hyperplane section. Now $H_X$ is a
smooth, non degenerate curve of genus $p$ and degree $2p$ in $\bold P^p$, with $p = 5$. The projective normality of such a model of a
genus $p$ curve is proved in [GL]. To prove that $X$ is rational observe that 
$$
K_X \sim H_X - F
$$
by 5.4, hence  $mH_XK_X = -2m$ and $P_m(X) = 0$ for $m \geq 1$. Moreover the exact sequence
$$
0 \to \Cal O_X \to \Cal O_X(H_X) \to \Cal O_{H_X}(H_X) \to 0
$$
implies that $X$ is regular. Indeed the restriction map $H^0(\Cal O_X(H_X)) \to H^0(\Cal O_{H_X}(H_X))$ is
surjective because $h^0(\Cal O_{H_X}(1)) = 6$ and $X, H_X$ are not degenerate. Passing to the long exact sequence we get
$$
0 \to H^1(\Cal O_X) \to H^1(\Cal O_X(H_X)) \to \dots
$$
On the other hand we compute $h^2(\Cal O_X(H_X)) = h^0(\Cal O_X(K_X-H_X)) = h^0(\Cal O_X(-E)) = 0$. Then Riemann-Roch yelds $h^1(\Cal
O_X(H_X)) = 0$. Hence $h^1(\Cal O_X) = 0$ and $X$ is regular. To show that the homogeneous ideal of
$X$ is generated by quadrics we consider the exact diagram
$$ 
\CD
0 @>>> {V \otimes H^0(\Cal I_{X \cup Y}(2))} @>>> {V \otimes (H^0(\Cal I_{X}(2) \oplus H^0(\Cal I_{Y}(2)))} @>>> {V
\otimes H^0(\Cal I_{F}(2))} @>>> 0
\\ @. @VVV @VV{u_X(2) \oplus u_Y(2)}V @VVV \\
0 @>>> {H^0(\Cal I_{X \cup Y}
(3))} @>>> {H^0(\Cal I_{X}(3)) \oplus H^0(\Cal I_{Y}(2))} @>>> {H^0(\Cal I_{F}(3)}) @>>> 0 \\
\endCD
$$
where the vertical arrows are the natural multiplication maps and $V = H^0(\Cal O_{\bold P^6}(1))$. The construction of the diagram is
easy: the starting point is diagram is tensoring by $V$ the Mayer-Vietoris sequence
$$
0 \to \Cal I_{X \cup Y}(2) \to \Cal I_X(2) \oplus \Cal I_Y(2) \to \Cal I_F(2) \to 0.
$$
Since $X \cup Y$ is a complete intersection,  $h^1(\Cal I_{X \cup Y}(m)) = 0, \ m \geq 1$, hence the diagram is exact.
We already know from the proof of the previous proposition that $F$ is a not trigonal canonical curve. Hence the right
vertical arrow is surjective. The same is true for the left arrow, because $X \cup Y$ is complete intersection. Thus
the central arrow is surjective. But this is the direct sum of the multiplication maps $u_X(2)$ and $u_Y(2)$,therefore 
both $u_X(2)$ and $u_Y(2)$ are surjective. To complete the proof it suffices to show that the multiplication
$$
u_X(m): V \otimes H^0(\Cal I_X(m)) \to H^0(\Cal I_X(m+1))
$$
is surjective for $m \geq 3$. By Castelnuovo-Mumford theorem $u_X(m)$ is surjective for $m \geq 3$ if $h^i(\Cal I(3-i)) = 0$, $i > 0$.
This follows from the standard exact sequence
$$
0 \to \Cal I_X(3-i) \to \Cal O_{\bold P^6}(3-i) \to \Cal O_X(3-i) \to 0,
$$
using the projective normality of $X$ and the vanishing of $h^1(\Cal O_X(2))$ and of $h^2(\Cal O_X(1))$. \bigskip \noindent 
 \bigskip \noindent
Let
$$
f: X \to \bold P^4 \tag 5.6
$$
be the adjoint map defined by the linear system 
$$
\mid K_X + H_X \mid = \mid 2H_X - F \mid.
$$
Since $X$ is projectively normal, the linear system defining $f$ is cut on $X$ by the quadrics containing $F$. We know from Noether's
theorem that the ideal of $F$ is generated by quadrics, since $F$ is not trigonal. Hence $ \mid 2H_X- F\mid $ is base-point-free and
$f$ is a morphism. Let
$$
S = f(X), \tag 5.7
$$
it is easy to compute $(K_X+H_X)^2 = 4$ and $p_a(H_X+K_X) = 1$. This implies that $f$ is birational onto $S$ and that $S$ is a quartic 
Del Pezzo surface. Applying to $f$ Reider's theorem it follows that $f$ contracts exactly six lines to the distinct points 
$$
b_1 \dots b_6 \in S. \tag 5.8
$$
(5.9) \bf LEMMA \it One can assume that $b_1, \dots, b_6$ are general points on $S$ and that $S$ is a general quartic Del
Pezzo surface.
\rm \bigskip \noindent
PROOF In the Hilbert scheme of quartic Del Pezzo surfaces consider the open set $U$ parametrizing integral surfaces.  Let $u: \Cal S \to U$
be the universal surface and let 
$$
u_6: \Cal S^6 \to U
$$
be the six times fibre product of $u$. $U$ is irreducible and the fibre of $u_6$ at $S$ is $S^6$, therefore $\Cal S^6$ is irreducible. Let
$b = (b,_1, \dots, b_6)$, then  $(b,S)$ is a point of $\Cal S^6$. Note that $(b,S)$ defines the 6-dimensional linear system
$$
\mid \Cal I_{Z_b/S}(2) \mid
$$
where $Z_b =: \lbrace b_1 \dots b_6 \rbrace$. The associated map $f_b: S \to \bold P^6$ is just $f^{-1}$ and the image
$$
X_b \tag 5.10
$$
of $f_b$ is just $X$. For the point $(b,S)$ the surface $X_b$ is smooth, projectively normal and its homogeneous ideal is generated by
quadrics. All these properties are preserved on an open neighborhood of $(b,S)$, therefore they are satisfied by $X_{b'}$ for a general
$(b',S') \in \Cal S^6$. This implies the statement. \bigskip \noindent Next we recall that there exists a blowing down
$\sigma: S \to \bold P^2 $ of 5 disjoint lines $L_1, \dots, L_5$ of $S$. By
the lemma we can assume that
$L_i \cap \lbrace b_1 \dots b_6 \rbrace$ is empty. Thus
$$
\tau =: \sigma \cdot f
$$ 
is the blowing up of 11 distinct points of $\bold P^2$, that is
$$
l_i = \sigma (L_i), \ i = 1 \dots 5 \ \text {and} \ e_j = \tau (E_j), \ j = 1 \dots 6, \tag 5.11
$$
where $E_j$ is the exceptional line contracted by $f$ to $b_j$. Let $P \in \mid \sigma^* \Cal O_{\bold P^2}(1) \mid $, note that 
$$
Pic X = \bold Z [P] \oplus \bold Z [L_1] \oplus \dots \bold Z [L_5] \oplus \bold Z [E_1] \oplus \dots \bold Z [E_6]. \tag 5.12
$$
It is easy to compute that
$$
\mid H_X \mid = \ \mid 6P - 2(L_1 + \dots + L_5) - (E_1 + \dots + E_6) \mid. \tag 5.13
$$
By the lemma we can assume that $l_1, \dots, l_5, e_1, \dots, e_6$ are sufficiently general, in particular that $l_1, l_2, e_1, e_2$ are
the base points of an irreducible pencil of conics. The strict transform on $X$ of a conic of this pencil will be denoted by
$$
R. \tag 5.14
$$
It is clear that $\mid R \mid$ is irreducible and base-point-free, a general $R$ is a smooth rational curve in $\bold P^6$ of degree
$6 = H_XR$. We assume from now on that $l_1 \dots l_5 e_1 \dots e_6$ are in general position in $\bold P^2$.
\bigskip \noindent
(5.15) \bf LEMMA \it (1) $R$ is non degenerate. \par \noindent
(2) Let $R' \subset R$ be a proper irreducible component, $R'$ is linearly normal.
\rm \bigskip \noindent
PROOF (1) Note that $H_X-R \sim 4P-2(L_3+L_4+L_5) - (L_1+L_2+E_3+E_4+E_5+E_6)$. Therefore $\mid H_X - R \mid$ is non empty if and only if 
there exists a quartic curve $Q \subset \bold P^2$ passing through $l_1, \dots, l_5, e_3, e_4, e_5, e_6$ and singular at $l_3,l_4,l_5$.
This does not happen if these points are sufficiently general in $\bold P^2$. 
\par \noindent
(2) $R'$ is either the strict transform of a line joining two of the points $l_1, l_2, e_1, e_2$ or the strict transform of a
smooth conic through $l_1,l_2,e_1,e_2$ and $o \in \lbrace l_3, l_4, l_5, e_3, e_4, e_5, e_6 \rbrace$. Let $o = e_i$,
$i = 3 \dots 6$. Then $R'$ is a smooth rational quintic and it suffices to show that $dim \mid H_X - R + E_6 \mid$  $= 0$.  This is
equivalent to say that there exists a unique plane quartic passing through $l_1, \dots, l_5, e_3, e_4, e_5$ and singular at $l_3, l_4,
l_5$: since $l_1, \dots, l_5, e_3, e_4, e_5$ are general this is true. We omit further details.
\bigskip
\noindent 
\bigskip \noindent
Finally a curve
$$
C \in \ \mid 2H_X - R \mid
$$
has arithmetic genus 8 and degree 14. $C$ is exactly the curve we are looking for: \bigskip \noindent 
(5.16) \bf THEOREM \it A smooth $C \in \mid 2H_X - R \mid$ belongs to $\Cal C_{14,8,6}$, moreover its homogeneous ideal is generated by
quadrics.
\rm \bigskip \noindent
PROOF A general $R$ as above is a smooth, non degenerate rational sextic curve in $\bold P^6$. Hence the homogeneous ideal
of $R$ is generated by quadrics. From this and the projective normality of $X$ it follows that 
$$
\mid 2H_X - R \mid
$$
is base-point-free. Then, by Bertini's theorem, a general $C \in \mid 2H-R \mid$ is smooth, and it is connected because $C^2 > 0$. 
Tensoring the standard exact sequence
$$
0 \to \Cal I_{X/\bold P^6}(2) \to \Cal I_{C/\bold P^6}(2) \to \Cal I_{C/X}(2) \to 0
$$
by $V = H^0(\Cal O_{\bold P^6}(1))$ and passing to the long exact sequence we obtain
$$
0 \to V \otimes H^0(\Cal I_{X/\bold P^6}(2)) \to V \otimes H^0(\Cal I_{C/\bold P^6}(2)) \to V \otimes H^0(\Cal I_{C/X}(2)) \to 0.
$$
The multiplication $l: V \otimes H^0(\Cal I_{X/\bold P^6}(2)) \to H^0(\Cal I_{X/\bold P^6}(3))$ is surjective because the ideal
of $X$ is generated by quadrics. On the other hand we have $\Cal I_{C/X}(2) \cong \Cal O_X(R)$. Thus, if the multiplication 
$$
r: V \otimes H^0(\Cal O_X(R)) \to H^0(\Cal O_X(H_X+R))
$$
has maximal rank, then the same is true for $\mu_C: V \otimes H^0(\Cal I_{C/\bold P^6}(2)) \to H^0(\Cal I_{C/\bold P^6}(3))$. 
From the exact sequence
$$
0 \to \Cal O_X(H_X) \to \Cal O_X(H_X+R) \to \Cal O_R(H_X+R) \to 0
$$
we obtain $h^0(\Cal O_X(H_X+R)) = 14$. Hence $r$ has maximal rank iff $r$ is an isomorphism. Now $\mid R \mid$ is a base-point-free
pencil, in particular it has no fixed divisor. Then, applying the base-point-free pencil trick as proved for curves in [ACGH] p.126, 
it follows that  $Ker \ r$ $=$ $H^0(\Cal O_X(H_X-R))$. Since $R$ is non degenerate $Ker \ r = (0)$. Hence $r$ is an isomorphism and
$\mu_C$ is surjective. Using this fact, and Castelnuovo-Mumford theorem as in 5.5, one deduces that the ideal of $C$ is 
generated by quadrics if $h^i(\Cal I_{C/\bold P^6}(3-i)) = 0$ for $i > 0$. This follows from the long exact sequence of
$$
0 \to \Cal I_{C/\bold P^6}(3-i) \to \Cal O_{\bold P^6}(3-i) \to \Cal O_C(3-i) \to 0
$$
if $h^1(\Cal I_{C/\bold P^6}(2)) = h^1(\Cal O_C(1)) = 0$. Since $h^1(\Cal O_X(R))$ $= h^1(\Cal I_{X/\bold P^6}(2)) = 0$, we
already have $h^1(\Cal I_{C/\bold P^6}(2)) = 0$. To show that $\Cal O_C(1)$ is non special consider the long exact sequence of
$$
0 \to \Cal O_X(R-H_X) \to \Cal O_X(H_X) \to \Cal O_C(1) \to 0.
$$
We have $h^1(\Cal O_X(H_X)) = 0$ then it suffices to show $h^2(\Cal O_X(R-H_X)) = 0$. By Serre duality this is $h^0 (\Cal O_X(K_X+H_X-R))$.
In $Pic(X)$ we have $K_X+H_X-R \sim P - (L_3 + L_4 + L_5) + E_1 + E_2$. Since $l_3, l_4, l_5$ are not collinear points the latter divisor
is not linearly equivalent to an effective one. Hence
$h^2(\Cal O_X(R-H_X))$ $=$ $0$. Finally $C$ is projectively normal: indeed $C$ is non degenerate because $H_X(H_X - C) < 0$, moreover 
the non speciality of $\Cal O_C(1)$ implies that $C$ is linearly normal. Since $C$ is also 2-normal the projective normality of $C$
follows, ([ACGH] p. 140 D-5). In particular we have also shown that $C$ is in $\Cal C_{14,8,6}$.
\bigskip \noindent
\bf 6 The unirationality of $\Cal M_{14}$. \rm In order to show the unirationality of $\Cal M_{14}$ we consider our usual Grassmann
bundle
$$
u: \Cal G^{\sigma}_{14,8,6} \to \Cal C_{14,8,6}
$$
where we put $\sigma = (2,5)$. Then, from the formulae given in 2.10, we compute that
$$
(d,g,r) = (14,8,6) \ \Leftrightarrow \ (d',g',r) = (18,14,6).
$$
By theorem 5.16 the homogeneous ideal of a general $C \in \Cal C_{14,8,6}$ is generated by quadrics. Moreover the condition
$2g' - 2 > d' > f^{r-2}$ is satisfied and the Brill-Noether number $\rho(d',g',r)$ is $0$. Then, applying theorem 4.5 to this
case, it follows: 
\bigskip \noindent
(6.1) \bf THEOREM \it Both $\Cal M_{14}$ and $\Cal U^6_{18,14}$ are unirational. \rm \bigskip \noindent
Note that, via Serre duality, the main Brill-Noether locus $\Cal U^6_{18,14}$ is biregular to $\Cal U^1_{8,14}$.
\bigskip \noindent
\bf 7 The unirationality of $\Cal M_{12}$. \rm We put again $\sigma = (2,5)$ and apply a very similar argument.  \bigskip \noindent
(7.1) CLAIM \it $\Cal G^{\sigma}_{15,9,6}$ satisfies the key assumption. \rm \bigskip \noindent
We note that $(d,g,r) = (15,9,6) \Leftrightarrow (d',g',r) = (17,12,6)$. Under the claim the map
$$
\beta_{17,12,6} \cdot \gamma_{15,9,6}: \Cal G^{\sigma}_{15,9,6} \to \Cal W_{17,12,6}
$$
exists, by proposition 4.2 its image is open. Again the condition $2g' - 2 > d' > f^{r-2}$ holds in this case and $\rho(d',g',r) = 5$ is
positive. Via Serre duality $\Cal W^6_{17,12,6}$ is biregular to $\Cal W^0_{5,12}$. This is the universal 5-symmetric product  hence it is
irreducible. Then, applying theorem 4.5, it follows:
\bigskip \noindent
(7.2) \bf THEOREM \it Both $\Cal W^0_{5,12}$ and $\Cal M_{12}$ are unirational. \rm \bigskip \noindent
PROOF OF THE CLAIM Let
$$
X \cup Y
$$
be the reducible complete intersection of 4 quadrics considered in section 4. Keeping the assumptions and notations used there,  we
consider on $X$ the irreducible, base-point-free pencil of rational normal sextics $\mid R \mid$. This pencil contains the curve
$$
D_1 + E_6 \in \mid R \mid, \tag 7.3
$$
where $D_1$ is the strict transform on $X$ of the irreducible conic passing through the points $l_1, l_2, e_1, e_2, e_6$. We know
from lemma 5.15-(2) that $D_1$ is a smooth rational normal quintic spanning a hyperplane in $\bold P^6$. Now we consider the linear 
system
$$
\mid 2H_X - D_1 \mid
$$
of curves of genus $9$ and degree 15. Since $D_1$ is a non degenerate rational normal quintic, the sheaf $\Cal I_{D_1/\bold P^6}(2)$ is
globally generated. Hence the image of the natural restriction map
$$
\rho: \mid \Cal I_{D_1/\bold P^6}(2) \mid \to \mid 2H_X - D_1 \mid
$$
is base-point-free. Since $(2H_X-D_1)^2 > 0$ it follows that a general $C \in \mid 2H_X - D_1 \mid$ is a smooth, irreducible
curve. Moreover we have:  \bigskip \noindent
(7.5) \bf THEOREM \it $C$ is projectively normal and $\Cal O_{C}(1)$ is non special, so that $C \in \Cal C_{17,12,6}$. \rm \bigskip
\noindent PROOF Note that $C \sim C' + E_6$, where $C' \in \mid 2H_X - R \mid$. We have already studied $\mid 2H_X - R \mid$: we know 
from theorem 5.16 and its proof that $\mid 2H_X - R \mid$ is base-point-free ant that a general $C' \in \mid 2H_X - R \mid$ is a
projectively normal element of $\Cal C_{14,8,6}$. In particular we can assume that such a general $C'$ is transversal to $E_6$. We also
recall that $E_6$ is a line and that $Z =: E_6 \cap C'$ is supported on 2 points. The non speciality of $\Cal O_{C'\cup E_6}(1)$ follows
from the long exact sequence of
$$
0 \to \Cal O_{C' \cup E_6}(1) \to \Cal O_{C'}(1) \oplus \Cal O_{E_6}(1) \to \Cal O_Z (1) \to 0.
$$
In a completely analogous way, the vanishing of $h^1(\Cal I_{C' \cup E_6}(m))$, $m \geq 1$, follows from the long exact sequence of
$$
0 \to \Cal I_{C' \cup E_6/\bold P^6}(m) \to \Cal I_{C'}(m) \oplus \Cal I_{E_6/ \bold P^6}(m) \to \Cal I_{Z/\bold P^6}(m) \to 0.
$$
Then, by semicontinuity the same properties hold for a general $C \in \mid C' + E_6 \mid$. \bigskip \noindent
Now we fix a general, smooth $C \in \mid 2H_X - D_1 \mid$, then $C$ is transversal to $D_1$ and
$$
C \cup D_1 
$$
is a nodal quadratic section of $X$. $X$ is a scheme theoretic intersection of quadrics, hence the same property holds for
$C \cup D_1$. Then, applying proposition 3.2, there exists a nodal complete intersection of 5 quadrics
$$
C \cup D_1 \cup D_2 = Q_1 \cap  \dots \cap Q_5. \tag 7.6
$$
From formulae 2.10 we compute that $D_1 \cup D_2$ is a nodal curve of degree 17 and arithmetic genus 12. On the other hand the
surface $X$ is linked to a smooth sextic Del Pezzo $Y$ by a complete intersection of 4 quadrics, so it is not restrictive to assume
$$
Q_1 \cap \dots Q_4 = X \cup Y.
$$
But then $D_2$ is a quadratic section of $Y$, hence it is smoothable to a canonical curve of genus 7. Observe that $Sing \ D_2 \cap (C \cup
D_1)$ is empty because $C \cup D_1 \cup D_2$ is nodal. Moreover $B =: D_1 \cap D_2$ is a set of six linearly independent points on the
rational normal quintic $D_1$. Then, keeping $D_1$ fixed and moving $D_2$ in $\mid \Cal I_{B/Y}(D_1) \mid$, we can smooth $D_2$. This
shows that there exists a flat family
$$
D_1 \cup D_{2,t}, \ t \in T,
$$
such that $D_{2,t}$ is a smooth canonical curve for $t \in T-o$ and $D_{2,o}Ê= D_2$. Up to shrinking $T$ we can assume that $h^0(\Cal
I_{D_1 \cup D_{2,t}}(2))$ is constantly equal to 5. This follows by semicontinuity from the next lemma. \bigskip \noindent
(7.7) \bf LEMMA \it $h^1(\Cal I_{D_1 \cup D_2/ \bold P^6}(2)) = 0$ and $h^0(\Cal I_{D_1 \cup D_2/\bold P^6}(2)) = 5$. \rm \bigskip
\noindent PROOF $D_2$ is a quadratic section of $Y$, so we have the standard exact sequence of ideal sheaves
$$
0 \to \Cal I_{Y/\bold P^6}(2) \to \Cal I_{D_2/\bold P^6}(2) \to \Cal O_Y(D_1) \to 0.
$$
The associated long exact sequence yelds $h^1(\Cal I_{D_2/\bold P^6}(2)) = 0$. Then the long exact sequence of
$$
0 \to \Cal I_{D_1 \cup D_2/\bold P^6}(2) \to \Cal I_{D_1/\bold P^6}(2) \oplus \Cal I_{D_2/\bold P^6}(2) \to \Cal I_{D_1 \cap D_2/\bold
P^6}(2) \to 0
$$
yelds $h^0(\Cal I_{D_1 \cup D_{2,t}/\bold P^6}(2)) = 5$. \bigskip \noindent
Since $h^0(\Cal I_{D_1 \cup D_{2,t}/\bold P^6}(2)) = 5$, the above complete intersection $C \cup D_1 \cup D_2$ deforms in a flat
family of complete intersections of 5 quadrics:
$$
C_t \cup D_1 \cup D_{2,t} = Q_{1,t} \cap \dots \cap Q_{5,t}, \ t \in T.
$$
Since $C \in \Cal C_{15,9,6}$, a general $C_t$ belongs to $\Cal C_{15,9,6}$. Since a general $D_t$ is smooth, we conclude that it is
not restrictive to assume that $D_2$ is a smooth canonical curve of genus 7. \bigskip \noindent 
(7.9) \bf LEMMA \it $h^1(T_{\bold P^6} \otimes \Cal O_{D_1 \cup D_2}) = 0$ so that $D_1 \cup D_2$ is smoothable. \rm \bigskip \noindent
PROOF Consider the Mayer-Vietoris sequence 
$$
0 \to T_{\bold P^6} \otimes \Cal O_{D_1 \cup D_2} \to T_{\bold P^6} \otimes (\Cal O_{D_1} \oplus \Cal O_{D_2}) \to T_{\bold P^6} \otimes
\Cal O_{D_1
\cap D_2}
\to 0.
\tag 7.9
$$
As is well known $h^1(\Cal T_{\bold P^6} \otimes \Cal O_{D_i}) = 0$ for the curves we are considering. On the other hand $D_1 \cap D_2$
is a set of 6 linearly independent points and $D_1$ is a non degenerate rational normal curve. Hence the restriction map $H^0(T_{\bold
P^6} \otimes \Cal O_{D_1}) \to T_{\bold P^6} \otimes \Cal O_{D_1 \cap D_2}$ is surjective. These remarks, and the long exact sequence of
7.9, imply $h^1(T_{\bold P^6} \otimes \Cal O_{D_1 \cup D_2}) = 0$. Then $D_1
\cup D_2$ is smoothable, (cfr. [HH]	2.1)
\bigskip \noindent
We conclude the proof of our claim: let $\lbrace D_t, t \in T \rbrace$ be a flat family of curves in $\bold P^6$ such that $D_t$ is smooth
for $t \in T-o$ and $D_o = D_1 \cup D_2$. Let $V_t = H^0(\Cal I_{D_t/\bold P^6}(2))$. By lemma 7.8  we can assume $dim \ V_t = 5$ 
for each $t \in T$. The scheme defined by $V_o$ is $C \cup D_1 \cup D_2$. Hence $V_t$ defines a complete intersection of quadrics $C_t \cup
D_t$, with $C_t \in \Cal C_{15,9,6}$ and the claim follows.

\rm \bigskip 
\noindent
\bf 8 The unirationality of $\Cal M_{11}$. \rm In this case we shift to curves in $\bold P^4$, but the arguments are the
same. We definitely assume $\sigma = (3,3)$ and consider the Grassmann bundle
$$
u: \Cal G^{\sigma}_{13,9,4} \to \Cal C_{13,9,4}.
$$
For $\sigma = (3,3)$ it turns out that
$$
(d,g,r) = (13,9,4) \Longleftrightarrow (d',g',r) = (14,11,4).
$$
(8.4) CLAIM \it $\Cal G^{\sigma}_{13,9,4}$ satisfies the key assumption. \rm \bigskip \noindent
This will be shown in theorem 9.10 of the next section. Under the claim the map 
$$
\gamma_{13,9,4}: \Cal G^{\sigma}_{14,11,4} \to \Cal W^4_{14,11}
$$
exists. By proposition 4.4 the image of $\beta_{14,11,4} \cdot \gamma_{13,9,4}: \Cal G^{\sigma}_{13,9,4} \to \Cal W^2_{14,11}$
is open. Serre duality yelds a biregular map between $\Cal W^4_{14,11}$ and the universal 6-symmetric product $\Cal W^0_{6,11}$.
Hence $\Cal W^4_{14,11}$ irreducible. Applying theorem 4.5 it follows: 
\bigskip \noindent
(8.5) \bf THEOREM \it Both $\Cal W^0_{6,11}$ and $\Cal M_{11}$ are unirational. \rm \bigskip \noindent
\bf 9. Curves of degree 12 and genus 8 in $\bold P^4$. \rm In this section we prove our previous claim 8.4. Preliminarily
we also show that a general curve of degree 12 and genus 8 in $\bold P^4$ is a scheme theoretic intersection of cubics. This will be
used in the next section. \par \noindent
Instead of our favourite rational surface of degree 10 in $\bold P^6$, we use now a smooth septic 
$$
X \subset \bold P^4
$$
having sectional genus 5 and birational to a K3 surface. This surface is very well known and its properties are described,
(cfr. [DES]). We preliminarily recall some of them: \bigskip \noindent
(9.1) \bf PROPOSITION \it $X$ is projectively normal and its ideal is generated by 3 cubic forms. \rm 
\bigskip \noindent
The geometric construction of $X$ is also well known, (cfr. [Ba]): \bigskip \noindent
(9.2) \bf PROPOSITION \it Let $X' \subset \bold P^5$ be any smooth complete intersection of 3 quadrics and let $e \in X'$
be a point not on a line of $X'$. Then the image of $X'$ under the linear projection of center $e$ is a smooth
surface $X$ as above. \rm \bigskip \noindent
Thus $X$ is defined by the blowing up $\sigma: X \to X'$ at $e$ and $K_X = \sigma^{-1}(e)$, we have
$$
Pic \ X = \sigma^* Pic \ X' \oplus \bold Z[K_X].
$$
To construct some curves of genus 8 and 9 we choose a suitable $X'$ with Picard number two: 
\bigskip \noindent
(9.3) \bf PROPOSITION \it There exists a smooth complete intersection of 3 quadrics $X' \subset \bold P^5$ such that
$$
Pic \ X' = \bold Z[L'] \oplus \bold Z[H'],
$$
where $H'$ is a hyperplane section of $X'$ and $L'$ is a very ample curve of degree 10 and genus 3. Moreover any
effective divisor on $X'$ is very ample.
\rm \bigskip \noindent 
PROOF The existence of a K3 surface $X'$ with Picard lattice as above is a standard consequence of the surjectivity  of the periods map 
for K3 surfaces. Such a lattice does not contain non zero vectors $v$ such that $v^2 = 0, -2, 2$. Indeed let $v = xH' + yL'$, then $v^2$ 
$=$ $4(2x^2 + y^2 + 5xy)$ $\neq 2, -2, 0$. Let $D$ be any effective divisor, then $D^2 \geq 4$ and $dim \ \mid D \mid$ $\geq 3$. Let $F$ 
be a fixed irreducible component of $\mid D \mid$, then $dim \mid F \mid = 0$ and hence $F^2 < 0$: a contradiction. Since $X'$ is a K3
surface, then $\mid D \mid$ is base-point-free and irreducible, moreover $\mid D \mid$ is very ample unless $D^2 = 2$ or there exists a
curve $F$ such that $DF \leq 2$ and $F^2 \in \lbrace 0, -2 \rbrace$. Hence $D$ is very ample. Up to changing their sign, we can assume that
both the generators $H'$ and $L'$ of $Pic \ X'$ are effective. In particular we can assume that $X'$ is embedded in $\bold P^5$ by $H'$.
Then either $X'$ is a complete intersection of 3 quadrics or contains a pencil $\mid F \mid$ of plane cubics. Since $F^2 = 0$ the latter
case is excluded.
\bigskip \noindent
From now on we assume that $X'$ is a K3 surface as in the previous statement. On $X'$ we have the very ample linear system 
$$
\mid 3H' - L' \mid. \tag 9.4
$$
of curves of degree 14 and genus 9. Let $e \in X'$ be a general point, due to the very ampleness of the linear systems we are
considering we can assume that:
\bigskip \noindent 
(1) There exists $A'_e \in \mid 3H' - L' \mid$ having an ordinary node at $e$ and no other singular point. \par
\noindent  (2) There exists $B'_e \in \mid L' \mid$ having an ordinary node at $e$ and no other singular point. \par \noindent  (3) The
linear systems $\mid L' - e \mid$ and $\mid 3H' - L' - e \mid$ have a unique, simple base point at $e$. \bigskip \noindent Finally let
$$
\pi: X' \to X \subset \bold P^4 \tag 9.5
$$
be the projection of center $e$, $\pi$ is the inverse of the blow up $\sigma: X \to X'$ at $e$. Let
$$
A \subset X \subset \bold P^4 \tag 9.6
$$
be the strict transform of $A'_e$ by $\sigma$. Then $A$ is a smooth, irreducible curve of genus 8 and degree 12. Let
$H =: \sigma^*H'$, $L =: \sigma^*L'$ then
$$
A \in \mid 3H - L - 2K_X \mid.
$$
Unfortunately $\Cal O_A(1)$ is special: this happens to every curve in $X$, since $X$ is regular and $h^1(\Cal O_X(1)) = 1$.
Nevertheless we can use $A$ to show the following \bigskip \noindent
(9.7) \bf THEOREM \it A general $C \in \Cal C_{12,4,8}$ is a scheme theoretic intersection of cubics. \rm \bigskip \noindent
PROOF It is easy to see that the Hilbert scheme $\Cal H_{12,8,4}$ is irreducible, in particular $\Cal C_{12,8,4}$ is
dense in $\Cal H_{12,8,4}$. Hence there exists a flat family $\lbrace C_t, t \in T \rbrace$ such that $C_t \in \Cal C_{12,8,4}$
if $t \neq o$ and $A = C_o$. Then, to prove the theorem, it suffices to show that: (1) $h^0(\Cal I_{C_t/\bold P^4}(3))$ is
constant on $T$, (2) $A$ is scheme theoretic intersection of cubics. To show (1) it suffices to show that $A$ is 3-normal.
This follows from the long exact sequence of
$$
0 \to \Cal I_{X/\bold P^4}(3) \to \Cal I_{A/\bold P^4}(3) \to \Cal O_X(3H-A) \to 0 \tag 9.8
$$
observing that $h^1(\Cal O_X(3H-A)) = 0$ and that $X$ is 3-normal. To show (2) recall that the ideal of $X$ is generated by cubics.
Then, to prove that
$A$ is a scheme theoretic intersection  of cubics, it suffices to show that $\mid 3H - A \mid $ is base-point-free. This follows 
because $\mid 3H - A \mid$ is the strict transform on $X$ of $\mid 3H' - A' - e \mid$, whose unique base point is $e$.
\bigskip \noindent
Now we turn to curves of degree 13 and genus 9: the linear system $\mid 3H - L - K_X \mid$ is just the strict transform by $\sigma$ of
$\mid 3H' - L' - e \mid$, hence a general
$$
C_o \in \mid 3H - L - K_X \mid \tag 9.9
$$
is smooth, irreducible of degree 13 and genus 9. Let
$$
B \subset X \subset \bold P^4
$$
be the strict transform of $B'_e$, $B$ is a smooth octic of genus 2. We can assume that $C_o \cup B$ is nodal, notice also that
$C_o + B$ is a cubic section of $X$. Since the ideal of $X$ is generated by 3 cubics, the ideal of $C_o \cup B$ is 
generated by 4. Then, by 3.2, there exists a nodal complete intersection 
$$
F_1 \cap F_2 \cap F_3 = C_o \cup B \cup B_1
$$
where $F_1, F_2, F_3$ are cubics. By 3.3 we can also choose $F_1, F_2$ so that
$$
F_1 \cap F_2 = X \cup Y
$$
where $Y$ is smooth. Then $Y$ is a quadric and $B_1$ is a smooth curve of type (3,3) on it. Using $C_o$ and the previous remarks we can
finally show that: \bigskip \noindent (9.10) \bf THEOREM \it A general $C \in \Cal C_{13,9,4}$ is linked to a smooth, irreducible curve by
a complete intersection of 3 cubics. In particular $\Cal G^{\sigma}_{13,9,4}$ satisfies the key condition. \rm \bigskip \noindent
PROOF It is easy to see that $\Cal H_{13,9,4}$ is irreducible. Therefore there exists an irreducible flat family  $\lbrace C_t, t \in T
\rbrace$ such that $C_t \in \Cal C_{13,9,4}$ if $t \neq o$  and $C_t = C_o$ if $t = o$. $C_o$ is 3-normal: the proof
is exactly the same used for the curve $A$ in the  proof of theorem 9.7. Let $W_t = H^0(\Cal I_{C_t/\bold P^4}(3))$, then
$W_t$ has constant dimension 4. Let $V_o \subset W_o$ be the subspace defining the complete intersection $C_o \cup B \cup B_1$, then we 
can move $V_o$ in a family $\lbrace V_t, t \in T \rbrace$ of 3-dimensional subspaces $V_t \subset W_t$. We can assume that $V_t$ defines
a nodal complete intersection of 3 cubics
$$
C_t \cup D_t
$$
and that $D_t$ is nodal, non degenerate, of degree 14 and arithmetic genus 11. Let $\Cal H$ be the complete  Hilbert
scheme of $D_t$, clearly there exists a rational map
$$
\gamma: \Cal G^{\sigma}_{13,9,4} \to \Cal H
$$
sending a general $(C,V) \in \Cal G^{\sigma}_{12,8,4}$ to $D$, $C \cup D$ being the complete intersection defined by $V$. However we only
know that $D_o = \gamma(C_o,V_o)$ is the nodal union of two smooth curves, so any $D$ in the image of $\gamma$ could be singular.
To complete the proof we show that this does not happen: \par \noindent We recall that, by 4.1, each $D$ in the image of
$\gamma$ is linearly normal and satisfies  $h^0(\Cal I_{D}(3)) = 3$. Let $D = \gamma(C,V)$, the latter property implies that $V$ $=$
$H^0(\Cal I_{D/\bold P^4}(3))$ and hence that $\gamma$ is birational onto its image. So we can compare dimensions. \par
\noindent  Let $D$ be general in the image of $\gamma$. $D$ is a flat deformation of $D_o$. $D_o$ has 6 nodes, moreover $D_o = B_o \cup
B_1$, where $B_o$, $B_1$ are smooth, ireducible and $B_o$ is not degenerate. Assume that $D$ is not smooth. Then, since $D$ is general, we
have the following cases: (1) $D$ is the nodal union of two smooth, irreducible curves, one of them not degenerate. $D$ has at
most 6 nodes, (2) $D$ is nodal, irreducible with at most 6 nodes. We discuss separately the two cases. \par
\noindent  (1) Let $f: \Cal D \to T$ be a flat family such that $\Cal D_t$ is general in the image of $\gamma$ and $\Cal D_o = D_o$. We can
assume that $T$ is smooth, irreducible and that each $\Cal D_t$ satisfies the condition in (1). Let $\Cal D$ be irreducible, then the two
irreducible components of a general $\Cal D_t$ must have the same degree and genus. This implies  that the degree is 7 and the genus is 3.
Let $\Cal F \subset \Cal H$ be the family of all 6-nodal curves $D_1 \cup D_2$ such that $D_i$ is a smooth septic of genus 3: we have $dim
\Cal F = 48$ $< dim \Cal G^{\sigma} = 60$. Hence the image of $\gamma$ is not in $\Cal F$. Assume now that $\Cal D$ is reducible, it is
easy to deduce that then a general $\Cal D_t$ is like $D_o$ i.e. it has 6 nodes and it is the union of a smooth canonical curve of
genus 4 and of a smooth octic of genus 3. Again this family of reducible curves has dimension $< dim \Cal G^{\sigma}_{13,9,4}$.
\par \noindent 
(2) We can consider an analogous family $f: \Cal D \to T$. In this case $\Cal D_t$ is irreducible, nodal with
$\nu \leq 6$ nodes if $t \neq o$ and $\Cal D_o = D_o$. Moreover $\Cal D_t$ is non degenerate and linearly normal. Let $\Cal F_{\nu}$ be the
corresponding family of irreducible, nodal curves of genus $11 - \nu$ and degree 14. It suffices to compute that $dim \Cal F_{\nu}$
$<$ $60$. This is a not difficult exercise.
\bigskip \noindent

\bigskip \noindent
\bf 10. The unirationality of $\Cal M_{13}$ \rm Let $\sigma = (3,3)$, continuing in the same vein we first point out that 
$\Cal G^{\sigma}_{12,14,4}$ satisfies the key condition as follows from theorem 9.7. Since
$$
(d,g,r) = (12,8,4) \Longleftrightarrow (d',g',r') = (15,14,4)
$$
a general $C \in \Cal C_{12,8,4}$ is linked to a smooth, irreducible, non degenerate $D$ of genus 14 by a complete intersection
of three cubics:
$$
C \cup D = F_1 \cap F_2 \cap F_3.
$$
$D$ has genus 14 and not 13, however we will turn very soon to curves $D$ having exactly one node. We already know that  
$$
C \cup D \subset S, \tag 10.1
$$
where $S$ is a smooth complete intersection of two cubics. Let $H$ be a hyperplane section of $S$: we know as well that, since $C$ is a
scheme theoretic intersection of cubics, then the linear system
$$
\mid 3H - C \mid
$$
is base-point-free. From $h^0(\Cal I_{C/\bold P^4}(3)) = 6$ it follows $dim \mid 3H - C \mid = 3$. Notice also that
$(3H - C)^2 = 11$. This implies that, for every $D \in \mid 3H - C \mid$, the linear series $\mid \Cal O_D(D) \mid$ is a complete
$g^2_{11}$ with no base points. Since the degree is 11, the associated morphism
$$
f_D: D \to \bold P^2
$$
is birational onto its image. We will say that a curve is \it uninodal \rm if it is nodal with a unique node. 
\rm \bigskip \noindent
From now on we assume that \it $D \in \mid 3H - C \mid$ is uninodal, \rm moreover we will denote as
$$
\Cal D \tag 10.3
$$
the family of all uninodal $D \in \Cal H_{15,14,4}$ such that there exists a complete intersection $C \cup D$ as above. 
Let $\nu: D' \to D$ be the normalization and let $n = \nu^* Sing \ D$, the morphism $f_D \cdot \nu: D' \in \bold P^2$ is 
defined by the line bundle
$$
L_{D'} =: \nu^* \Cal O_D(D).
$$
Since $\mid L_{D'} \mid$ is base-point-free of degree 11, the curve $f_D(D)$ is an element of the quasi-projective variety of all 
reduced, irreducible  linearly normal plane curves of degree 11 and genus 13. Due to the fundamental results of the Brill-Noether 
theory, such a variety is irreducible and a non empty open set of it parametrizes nodal curves, (cfr. [HM1] p.40). Moreover its image
$$
U^2_{11,13} \subset \Cal W^2_{11,13} \tag 10.4
$$
dominates the moduli space $\Cal M_{13}$: this follows because the Brill-Noether number $\rho(11,13,2)$ is $\geq 0$. 
The image of the element $f_D(D)$ is just the moduli point of the pair $(D',L_{D'})$ in $\Cal W^2_{11,13}$. Since
$$
L_{D'} = \omega_{D'}(n) \otimes \nu^* \Cal O_D(-1),
$$
we can define a morphism
$$
\phi: \Cal D \to \Cal U^2_{11,13} \tag 10.5
$$
sending a uninodal $D$ to the moduli point of $(D',L_{D'})$. \bigskip \noindent
(10.6) \bf LEMMA \it $\phi$ is dominant. \rm \bigskip \noindent
PROOF We fix an irreducible flat family of pairs $\lbrace (D_t, L_t), t \in T \rbrace $ such that: (1) $D_t$ is a smooth, irreducible curve
of genus
$13$ and $L_t \in Pic^{11}(D)$ is globally generated with $h^0(L_t) = 3$, (2) $T$ dominates $\Cal U^2_{11,13}$ via the natural map, (3)
for $t = o$ $(D_o,L_o) = (D',L_{D'})$. Up to a finite base change we can also assume that: (4) on each $D_t$ there exists a rationally
determined effective divisor $n_t$ which is contracted to a point by the morphism $f_t: D_t \to \bold P^2$ defined by $L_t$, (5)
$n_o = n$. Now we consider the other family of pairs
$$
(D_t, \omega_{D_t}(n_t) \otimes L_t^{-1})
$$
and the corresponding family of associated maps $h_t: D_t \to \bold P^4$. To show that $\phi$ is dominant it suffices to show that
$h_t(D_t) \in \Cal D$ for $t$ general. Since $D = h_o(D_o)$ is uninodal, the general $h_t(D_t)$ is uninodal. Then it
is immediate to compute that $h^0(\Cal I_{h_t(D_t)/\bold P^4}(3)) \geq 3$. On the other hand we have $h^0(\Cal I_{D/\bold P^4}(3)) = 3$,
as follows applying to $C \cup D$ proposition 4.1. Thus, by semicontinuity, the same property holds for a general $D_t$. Finally let
$V_t$
$=$
$H^0(\Cal I_{h_t(D_t)/\bold P^4}(3))$: the scheme defined by $V_o$ is the curve $C \cup D$. Hence the scheme defined by a general  $V_t$
is a nodal curve $C_t \cup D_t$, with $C_t \in \Cal C_{12,8,4}$. Then $h_t(D_t) \in \Cal D$ and $\phi$ is dominant. 
\rm \bigskip \noindent  
\bigskip \noindent
(10.7) \bf DEFINITION \it A point $(C,V) \in \Cal G^{\sigma}_{12,8,4}$ is uninodal if $D = \gamma_{12,8,4}(C,V)$ is uninodal. 
\par \noindent The family of all uninodal points $(C,V)$ will be denoted by $\Cal N$. \rm \bigskip \noindent
(10.8) \bf LEMMA \it $\Cal N$ and $\Cal D$ are unirational. \rm \bigskip \noindent
PROOF It is clear that $\gamma_{12,8,4}(\Cal N) = \Cal D$. Hence it will be sufficient to show that $\Cal N$ is unirational. \par \noindent
Fix a general pair $(C,x) \in \Cal C_{12,8,4} \times \bold P^4$, then consider the family $F(C,x)$ of all uninodal points $(C,V)$ 
such that $x = Sing \ D$ and $\gamma_{12,8,4}(C,V) = D$. It is easy to see that $F(C,x)$ is biregular to an open subset of the 
Grassmannian
$$
G(2,I_x/I_{2,x}),
$$
where $I_x =: \lbrace f \in I \ / \ f \in m_x \rbrace$, $  \ I_{2,x} =: \lbrace f \in I_x \ / \ f \in m^2_x \rbrace \ $ and $I =:
H^0(\Cal I_{C/\bold P^4}(3))$. Moreover $F(C,x)$ is the fibre at the point $(C,x)$ of the morphism
$$
\pi: \Cal N \to \Cal C_{12,8,4} \times \bold P^4
$$
sending $(C,V)$ to $(C,x)$, with $x = Sing \ D$ and $D = \gamma_{12,8,4}(C,V)$. We prove that $\pi(\Cal N)$ is dense: let $\pi(C,V)$
$= (C,x)$, then there exists a nodal complete intersection of 3 cubics $C \cup D$ such that $Sing \ D = x$. This condition is open
on $(C,V)$ and on $x$, hence it holds on open neighborhoods $U_C$ of $C$ and $U_x$ of $x$. Therefore $U_C \times U_x \subset \pi(\Cal N)$
and $\pi(\Cal N)$ is dense. On a non empty open set $A$ of $\pi(\Cal N)$ the space $I_x/I_{2,x}$ has constant dimension. It is standard to
construct on $A$ a vector bundle $\Cal I$ having fibre $I_x/I_{2,x}$ at the point $(C,x)$. On the other hand the fibre
$F(C,x)$ of  $\pi$ is open in $G(2,I_x,I_{2,x})$. Hence $\Cal N$ is biregular to an open set of the Grassmann bundle
$G(2,\Cal Q)$. This is birational to $\Cal C_{12,8,4} \times \bold P^4 \times G(2,4)$, therefore it is unirational.
\bigskip \noindent As a straightforward consequence of the lemma we have: \bigskip \noindent  (10.9) \bf THEOREM \it $\Cal M_{13}$ is
unirational as well as the Severi variety of nodal plane curves of degree 11 and genus 13. \rm \bigskip \noindent
\bigskip \noindent
\bf 11. References \rm \par \noindent
 [ACGH] E. Arbarello, M. Cornalba, P. Griffiths, J. Harris {\it
Geometry of Algebraic Curves I} Springer-Verlag, Berlin (1984), 1-386 \par \noindent
[CR1] M.C. Chang, Z. Ran  {\it Unirationality of the moduli space of curves of genus 11, 13
(and 12) ,} Invent. Math. {\bf 76 } (1984) 41-54,
\par \noindent 
[EH] D. Eisenbud, J. Harris  {\it The Kodaira dimension of the moduli space of curves of genus
$\leq$ 23 ,} Invent. Math. {\bf 74 } (1983), 359-387 \par \noindent
[Ba] I. Bauer {\it Inner projections of algebraic surfaces: a finiteness result}, Crelle J. reine Angew. Math. {\bf 480} (1993)
\par \noindent  
[DES] W. Decker, L. Ein, F.O. Schreyer {Construction of surfaces in $\bold P^4$}, J. Alg. Geom. {\bf 2} (1995) 
\par \noindent
[Fa] G. Farkas  
{\it The geometry of the moduli space of curves of genus 23 } Math. Annalen {\bf 318} (2000) 43-65, \par 
\noindent
[FP] G. Farkas, M. Popa {\it Effective divisors on $M_g$ and a counterexample to the Slope
Conjecture } preprint (2002)
\par \noindent
[GL] M. Green, R. Lazarsfeld  
{\it On the projective normality of complete linear series on an algebraic curve,} \ Invent. \ Math. {\bf 83 } (1986) 
73-90\par
\noindent
[HH] R. Hartshorne, A. Hirschowiz  
{\it Smoothing Algebraic Space Curves,} in Algebraic Geometry, Sitges 1983 (E. Casas-Alvero, G.E. Welters, S.
Xambo-Descamps eds.)  L.N.M. {\bf 1124 } (1985), 98-131 \par
\noindent  
[HM1] J. Harris, I. Morrison 
{\it Moduli of Curves ,} Springer-Verlag (1991), 1-366 \par \noindent 
[HM2]  
{\it Slopes of effective divisors on the moduli space of curves ,}  Invent. Math. {\bf 99}
(1990), 321-335 \par
\noindent  
[HM] J. Harris, D. Mumford  {\it On the Kodaira dimension of the moduli space of curves ,} 
Invent. Math.{\bf 67 } (1982), 23-86 \par \noindent
[M] S, Mukai 
{\it Curves, K3 surfaces and Fano manifolds of genus $\leq 10$ ,} in 'Algebraic Geometry in honor of Masayoshi
Nagata, (H. Hijikata, H. Hironaka eds.) Kikokuniya, Tokio   (1988), 367-377 \par \noindent  
[MM] S. Mori, S. Mukai 
{\it The uniruledness of the moduli space of curves of genus 11 ,} in Algebraic Geometry, Proceedings Tokio/Kyoto 1982
L.N.M {\bf 1016 } (M.Raynaud, T. Shioda eds.) (1983), 334-353\par
\noindent  
[S1] E. Sernesi  {\it L'unirazionalit\'a della variet\'a dei moduli delle curve di genere 12,}
Ann. Sc. Norm. Sup. Pisa {\bf 8} (1981), 405-439
\par \noindent
[S2] E. Sernesi  {\it On the existence of certain families of curves,} Invent. Math. {\bf 75 }
(1984), 25-57
\par \noindent 
 [Se] F. Severi  {\it Vorlesungen ueber Algebraische Geometrie ,}  (E. Loeffler uebersetzung),
Teubner, Leipzig {\bf } (1921),\par \noindent
[ST] F.O. Schreyer, F. Tonoli  {\it Needles in a haystack: special varieties via
small fields,} in Mathematical computations with Macaulay
2,  (D. Eisenbud, D. Grayson, M. Stillman, B. Sturmfels eds.),  Springer-Verlag, Berlin (2002). 
\bigskip \noindent 
ADDRESS A. Verra Dipartimento di Matematica, Largo S. Leonardo Murialdo 1, 00146 ROME (Italy)
( verra\@mat.uniroma3.it)

\enddocument